\documentclass[a4paper,12pt]{article}

\usepackage[french, ngerman, italian, english]{babel}
\usepackage{amsmath,amsthm, amssymb}
\usepackage{setspace}
\usepackage{anysize}
\usepackage{url}
\usepackage{graphicx}
\usepackage[colorlinks=true]{hyperref}
\usepackage{mathptmx}
\usepackage{paralist}
\usepackage{verbatim}
\usepackage[utf8]{inputenc}

\theoremstyle{plain}
\newtheorem{thm}{\bf Theorem}[section]
\newtheorem{thmnonumber}{\bf Main Theorem}

\newtheorem{lemma}[thm]{\bf Lemma}
\newtheorem{corollary}[thm]{\bf Corollary}

\theoremstyle{definition}
\newtheorem{definition}{\bf Definition}[section]
\theoremstyle{remark}

\def \mm{{\mathfrak m}}
\def \aa{{\bf a}}
\def \bb{{\bf b}}
\def \cc{{\bf c}}
\def \NN{\mathbb N}
\def \AG{\bar{A}(G)}

\begin{document}

%
%\title{Unmixed Graphs that are Domains}
%\author{Bruno Benedetti \ and \ Matteo Varbaro}
%\address{Institut f\"{u}r Mathematik, MA 6-2, TU Berlin, Germany}
%\email{\url{benedetti@math.tu-berlin.de}}
%\address{Dipartimento di Matematica, Via Dodecaneso 35, I-16146, Universit\'a  degli Studi di Genova, Italy}
%\email{\url{varbaro@dima.unige.it}} 
%\subjclass{05C70, 13A30}
%\date{}
%\keywords{Minimal vertex covers, Symbolic powers, Unmixed graphs}
% \maketitle

\title{Unmixed Graphs that are Domains.}
\author{Bruno Benedetti\thanks{%
Supported by DFG via the Berlin Mathematical School.} \\
\footnotesize Institut f\"{u}r Mathematik, MA 6-2\\
\footnotesize TU Berlin, Germany\\
\footnotesize \url{benedetti@math.tu-berlin.de} 
\and \setcounter{footnote}{3} Matteo Varbaro \\
\footnotesize Dipartimento di Matematica\\
\footnotesize Univ. degli Studi di Genova, Italy\\
\footnotesize \url{varbaro@dima.unige.it}}
\date{{\small \today}} 
\maketitle

\begin{abstract}
\noindent
 We extend a theorem of Villareal on bipartite graphs to the class of all graphs. On the way to this result, we study the \emph{basic covers algebra} $\AG$ of an arbitrary graph $G$. We characterize with purely combinatorial methods the cases when: 1) $\AG$ is a domain, 2) $G$ is unmixed and $\AG$ is a domain. 
\end{abstract}

\section{Introduction and Notation}
Fix a graph on $n$ vertices and give each vertex a price; let the ``cost'' of an edge be the sum of the costs of its endpoints. A nonzero price distribution such that no edge is cheaper than $k$ euros is called a \textit{$k$-cover}. 

A $k$-cover and a $k'$-cover of the same graph can be summed vertex-wise, yielding a $(k+k')$-cover; one says that a $k$-cover is \emph{basic} if it cannot be decomposed into the sum of a $k$-cover and a $0$-cover. Basic $1$-covers of a graph are also known as ``minimal vertex covers'' and have been extensively studied by graph theorists. 

A graph $G$ is called \textit{unmixed} if all its basic $1$-covers have the same number of ones. For example, a square is unmixed, a pentagon is unmixed, yet a hexagon is not unmixed. A graph $G$ is called a \textit{domain} if, for all $s, t \in \NN$, any $s$ basic $1$-covers and any $t$ basic $2$-covers always add up to a basic $(s+2t)$-cover. For example, the square is a domain, while the pentagon and the hexagon are not domains. 

This notation is motivated by the following algebraic interpretation (see Herzog et al. \cite{HHT}, or Benedetti et al. \cite{BCV} for details). Let $S$ be a polynomial ring of $n$ variables over some field and let $\mm$ be its irrelevant ideal. Let $I(G)$ be the ideal of $S$ generated by all the monomials $x_ix_j$ such that $\{i,j\}$ is an edge of $G$. The ideal $I(G)$ is called \emph{edge ideal} of $G$, and its Alexander dual $J(G)=\cap_{\{i,j\}}(x_i,x_j)$ is called \emph{cover ideal} of $G$. The \emph{symbolic fiber cone} of $J(G)$ is $\AG=R_s(J(G))/\mm R_s(J(G))$, where $R_s(J(G))=\oplus_{i\geq 0} J(G)^{(i)}$ is the symbolic Rees algebra of $J(G)$ and $J(G)^{(i)}$ is the $i$th symbolic power of $J(G)$. 

The following facts are easy to prove:
\begin{compactitem}[ -- ]
\item $G$ is unmixed if and only if $I(G)$ is unmixed;
\item $G$ is a domain if and only if $\AG$ is a domain;
\item if $G$ is a domain, then $\AG$ is a Cohen-Macaulay algebra (and a normal domain). 
\end{compactitem} 
\medskip
\noindent In the present paper, we introduce three entirely combinatorial properties, called ``square condition'' [SC], ``weak square condition'' [WSC], and ``matching square condition'' [MSC].
\begin{definition} We say that a graph $G$ satisfies
\begin{compactitem}[ -- ]
\item \textbf{SC}, if for each triple of consecutive edges $\{i,i'\}$, $\{i,j\}$ and $\{j,j'\}$ of $G$, one has that $i' \neq j'$ and $\{i',j'\}$ is also an edge of $G$.
\item \textbf{WSC}, if $G$ has at least one edge, and for every non-isolated vertex $i$ there exists an edge $\{i,j\}$ such that for all edges $\{i,i'\}$ and $\{j,j'\}$ of $G$, $\{i',j'\}$ is also an edge of $G$. (In particular, $i' \neq j'$).
\item \textbf{MSC}, if the graph $G^{red}$ obtained by deleting all the isolated vertices of $G$ is non-empty, and admits a perfect matching such that for each edge $\{i, j\}$ of the matching, and for all edges $\{i,i'\}$ and $\{j,j'\}$ of $G$, one has that $\{i',j'\}$ is also an edge of $G$. (In particular, $i' \neq j'$).
\end{compactitem}
\end{definition}

We will see in Lemma \ref{thm:SC} that the first property is satisfied only by bipartite complete graphs, by isolated vertices, and by disjoint unions of these two types of graphs. 

The second property, WSC, is a weakened version of the first one. It was studied in \cite{BCV}, where the authors proved that when $G$ is bipartite, $G$ satisfies WSC if and only if $G$ is a domain. We extend this result to non-bipartite graphs:

\begin{thmnonumber} [Theorem \ref{thm:WSC}]  \label{WSC}
A graph $G$ satisfies WSC if and only if $G$ is a domain.
\end{thmnonumber}

Finally, the third property was investigated by Villarreal \cite[Theorem 1.1]{V}, who proved that, when $G$ is bipartite, $G$ satisfies MSC if and only if $G$ is unmixed.  In the present paper, we extend Villarreal's theorem to the non-bipartite case, showing that

\begin{thmnonumber} [Theorem \ref{thm:MSC}] \label{MSC} 
A graph $G$ satisfies MSC if and only if $G$ is an unmixed domain.
\end{thmnonumber}

This implies Villarreal's result because in the bipartite case all unmixed graphs are domains. However, many graphs satisfy MSC without being bipartite (see Theorem \ref{thm:technique1}). From an algebraic point of view, Main Theorem \ref{MSC} characterizes the graphs $G$ for which every symbolic power of the cover ideal of $G$ is generated by monomials of the same degree.

%Notice also that the item (A) above is a stronger requirement than unmixedness: for %example, the complete graph $K_4$ is unmixed and has a matching, but it does not satisfy %(A) and in fact it is not even a domain.

We point out that the proof of Main Theorem \ref{MSC} is not an extension of Villarreal's proof. We follow a different approach, introducing the graph \textbf{$G^{0-1}$}, which is obtained from $G$ by removing the isolated vertices and then by removing all edges $\{i,j\}$ such that there exists a basic $1$-cover $\aa$ of $G$ for which $a_i + a_j = 2$. This graph $G^{0-1}$ always satisfies SC (see Lemma \ref{thm:SCG}). Furthermore, $G^{0-1}$ has no isolated points if and only if the original graph $G$ satisfied WSC (see Theorem \ref{thm:WSC}); finally, $G^{0-1}$ admits a matching if and only if $G$ satisfies MSC (see Theorem \ref{thm:MSC}).

For example, let $G$ be the graph on six vertices, given by the edges $\{1,2\}$, $\{2,3\}$, $\{3,4\}$, $\{1,4\}$, $\{2,5\}$, $\{4,5\}$ and $\{5,6\}$. This graph $G$ has three basic $1$-covers. As $G^{0-1}$ is the disjoint union of a $K_{2,2}$ and a $K_{1,1}$, we have that $G$ is an unmixed domain and satisfies MSC.

\section{Proofs of the main theorems.}

\begin{lemma} \label{thm:lemma1}
Let $\{i,j \}$ be an edge of a graph $G$. For any integer $d \geq 1$, the following are equivalent:
\begin{compactenum}[\rm (1)]
\item for any $k \in \{1, \ldots, d\}$, and for any basic $k$-cover $\aa$, $a_i + a_j = k$;
\item for any basic $1$-cover $\aa$, $a_i + a_j = 1$; 
\item if $\{i,i'\}$ is an edge of $G$ and $\{j,j'\}$ an edge of $G$, then $i' \neq j'$ and $\{i',j'\}$ is also an edge of $G$.
\end{compactenum}
\end{lemma}

\begin{proof}
(2) is a special case of (1). To see that (2) implies (3) we argue by contradiction. If $G$ contains a triangle $\{i,i'\}$, $\{i,j\}$, $\{j,i'\}$, we claim that there is a basic $1$-cover $\aa$ such that $a_i=a_j = 1$ and $a_{i'}=0$. In fact, define a $1$-cover $\bb$ by setting $b_{i'}=0$, and $b_{k}=1$ for all $k \neq i'$. In case $\bb$ is basic we are done; otherwise, $\bb$ breaks into the sum of a basic $1$-cover $\aa$ and some $0$-cover. This $\aa$ has still the property of yielding $0$ on $i'$, and thus $1$ on $i$ and $j$, so we are done. If instead $G$ contains four edges $\{i,i'\}$, $\{i,j\}$ and  $\{j,j'\}$, but not the fourth edge $\{i',j'\}$, we claim that there is a basic $1$-cover $\aa$ such that $a_i=a_j = 1$ and $a_{i'}=a_{j'}=0$. The proof is as before: First one defines a $1$-cover $\bb$ by setting $b_{i'}=b_{j'}=0$, and $b_{k}=1$ for all $k$ such that $j' \neq k \neq i'$; then one reduces $\bb$ to a basic $1$-cover. 

Finally, assume (1) is false: Then there is a basic $k$-cover $\aa$ such that $a_i + a_j > k$. For the cover to be basic, there must be a neighbour $i'$ of $i$ such that $a_{i'} + a_i = k$, and a neighbour $j'$ of $j$ such that $a_{j'} + a_j = k$. But then $a_{i'} + a_{j'} = 2k - a_i - a_j < k$, so $\{i', j'\}$ cannot be an edge of $G$: Hence, (3) is false, too. Thus (3) implies (1).
\end{proof}

In the proof of the next Lemma we use a convenient shortening: we say that a $k$-cover $\aa$ can be ``lopped at the vertex $i$'' if replacing $a_i$ with $a_{i}-1$ in the vector $\aa$ still yields a $k$-cover.

\begin{lemma} \label{thm:lemma2} Let $G$ be a graph. $G$ is a domain if and only if $G$ has at least one edge, and for each non-isolated vertex $i$ there exists a vertex $j$ adjacent to $i$ in $G$ such that:
\begin{compactitem}[ -- ]
\item for any basic $1$-cover $\aa$ one has $a_i + a_j = 1$, and
\item for any basic $2$-cover $\bb$ one has $b_i + b_j = 2$.
\end{compactitem}
\end{lemma}

\begin{proof} The fact that $G$ is a domain rules out the possibility that $G$ might be a disjoint union of points; so let us assume that $G$ has at least one edge. $G$ is \textit{not} a domain if and only if a non-basic $(s+2t)$-cover of $G$ can be written as the sum of $s$ basic $1$-covers and $t$ basic $2$-covers, if and only if there is a vertex $i$ such that a certain sum $\cc$ of $s$ basic $1$-covers and $t$ basic $2$-covers can be lopped at the vertex $i$ (and in particular, this $i$ cannot be isolated), if and only if there exists a non-isolated vertex~$i$ such that, for each edge $\{i,j\}$, either there exists a basic $1$-cover $a$ such that $a_i + a_j > 1$, or there exists a basic $2$-cover $b$ such that $b_i + b_j > 2$.
\end{proof}

\begin{lemma}\label{thm:SC}
Let $G$ be a connected graph.
$G$ satisfies SC if and only if $G$ is either a single point or a $K_{a,b}$, for some $b \geq a \geq 1$.
\end{lemma}

\begin{proof}
The fact that a $K_{a,b}$ satisfies SC is obvious. For the converse implication, first note that a graph $G$ satisfying SC cannot contain triangles; moreover, if $G$ contained a $(2d+1)$-cycle, by SC we could replace three edges of this cycle by a single edge, hence $G$ would contain a $(2d-1)$-cycle as well. By induction on $d$ we conclude that $G$ contains no odd cycle. So $G$ is bipartite: If $[n]=A \cup B$ is the bipartition of its vertices, we claim that any vertex in $A$ is adjacent to any vertex in $B$. In fact, if $G$ has no three consecutive edges, then $G$ is either a point or a $K_{1,b}$ (for some positive integer $b$) and there is nothing to prove. Otherwise, take $a \in A$ and $b \in B$: Since $G$ is connected, there is an (odd length) path from $a$ to $b$. By SC, the first three edges of such a path can be replaced by a single edge, yielding a path that is two steps shorter. Iterating the trick, we eventually find a path of length~$1$ (that is, an edge) from $a$ to $b$.
\end{proof}

\begin{definition}
Let $G$ be a graph with at least one edge. We denote by $\mathbf{G^{0-1}}$ the graph that has:
\begin{compactitem}[ -- ]
\item as vertices, the vertices of $G^{red}$;
\item as edges, the edges $\{i,j\}$ of $G$ such that for every basic $1$-cover $\aa$ of $G$ one has $a_i + a_j = 1$.
\end{compactitem}
\end{definition}

\begin{lemma} \label{thm:SCG}
Let $G$ be an arbitrary graph with at least one edge. Then $G^{0-1}$ satisfies SC.
\end{lemma}

\begin{proof}
Assume $\{h,i\}$, $\{i,j\}$ and $\{j,k\}$ are three consecutive edges of $G^{0-1}$. For any basic $1$-cover $\aa$, $a_h + a_i = 1$, $a_i + a_j= 1$ and $a_j + a_k = 1$. Summing up the three equations we obtain that $a_h + 2a_i + 2a_j + a_k = 3$, thus $a_h +a_k=1$: so all we need to prove is that $\{h,k\}$ is an edge of $G$. But by Lemma \ref{thm:lemma1}, this follows from the fact that for any basic $1$-cover $\aa$ one has $a_i + a_j= 1$. 
\end{proof}

\begin{thm} \label{thm:WSC} Let $G$ be a graph with at least one edge. Then the following are equivalent:  
\begin{compactenum}[\rm (1)]
\item $G$ satisfies WSC;
\item $G$ is a domain;
\item $\AG$ is a domain;
\item $G^{0-1}$ has no isolated points.
\end{compactenum}
\end{thm}

\begin{proof}
 The equivalence of (1) and (2) follows from combining Lemma \ref{thm:lemma1} and Lemma \ref{thm:lemma2}. The equivalence of (2) and (3) was explained in the Introduction. The equivalence of (1) and (4) is straightforward from Lemma \ref{thm:lemma1}. 
\end{proof}

\begin{lemma} \label{thm:vabbaslemma1}
Let $G$ be a domain. Let $H_1$, $\ldots$, $H_k$ be the connected components of $G^{0-1}$; let $A_i \cup B_i$ be the bipartition of the vertices of $H_i$, for $i=1, \ldots, k$. Then:
\begin{compactenum}[\rm (1)]
\item if $G$ contains a triangle, all three edges are not in $G^{0-1}$; in particular, two vertices of the same $B_i$ (or of the same $A_i$) are not adjacent in $G$;
\item if $G$ contains an edge between a vertex of $A_i$ and a vertex of $A_j$ then it contains also edges between any vertex of $A_i$ and any vertex of $A_j$;
\item if $G$ contains an edge between $A_i$ and $A_j$, then it contains no edge between $B_i$ and $B_j \,$;
\item if $G$ contains an edge between $A_i$ and $A_j$ and another edge between $B_i$ and $B_k$ 
then it contains an edge between $A_j$ and $B_k$; 
\item if $G$ contains an edge between $A_h$ and $A_i$ and another edge between $A_h$ and $A_j$ then $G$ contains no edge between $B_i$ and $B_j$.
\end{compactenum}
\end{lemma}

\begin{proof}
(1): Choose a vertex $v$ of the triangle, and take a basic $1$-cover that yields $0$ on $v$. This $1$-cover yields $1$ on the other two vertices, so the edge opposite to $v$ does not belong to $G^{0-1}$. The second part of the claim follows from the fact that each $H_i$ is complete bipartite: Were two adjacent vertices both in $A_i$ (or both in $B_i$), they would have a common neighbour in $G^{0-1}$, so there would be a triangle in $G$ with two edges in $G^{0-1}$, a contradiction.

(2): take a vertex $i$ of $A_i$ adjacent in $G$ to a point $j$ of $A_j$. Let $i'$ be any point of $A_i$ different from $i$. By contradiction, there is a vertex $j'$ of $A_j$ that is not adjacent to $i'$. Construct a basic $1$-cover $\cc$ that yields $0$ on $i'$, and $0$ on $A_j$ ($\cc$ is well defined because no two vertices of $A_j$ can be adjacent, by the previous item). Since $c_j = 0$, $c_i$ must be $1$; and since $c_{i'}=0$, $\cc$ yields $1$ on $B_i$; but then all edges $\{i, b\}$ with $b \in B_i$ are not in $G^{0-1}$, a contradiction.

(3): assume $G$ contains an edge $\{i, j\}$ between $A_i$ and $A_j$, and choose any vertex $x$ of $B_i$, and any vertex $y$ of $B_j$. Take a basic $1$-cover $\aa$ such that $a_i + a_j = 2$. Since $\{x,i\}$ and $\{ y,j\}$ are in $G^{0-1}$, $a_i + a_x = a_j + a_y = 1$; thus $a_x = a_y = 0$, which implies that there cannot be an edge in $G$ from $x$ to $y$. 

(4): fix a vertex $i$ of $A_i$ and use the WSC property (which $G$ satisfies by Theorem \ref{thm:WSC}): There exists a $x$ adjacent to $i$ such that for any edge $\{i, j\}$ and for any edge $\{x, y\}$ of $G$, $\{j, y\}$ is also an edge of $G$. By Lemma \ref{thm:lemma1}, $a_i + a_x = 1$ for each basic $1$-cover $\aa$; that is to say, $\{i,x\}$ is in $G^{0-1}$. This implies that $x$ is in $B_i$. So if $G$ contains an edge $\{i, j\}$, with $j \in A_j$, and an edge $\{x, y\}$, with $y$ in some $B_k$, then $G$ contains also the edge $\{j, y\}$ from $A_j$ to $B_k$.

(5): by contradiction, assume there is an edge between $B_i$ and $B_j$. By the previous item, since there is an edge between $A_h$ and $A_j$, there is also an edge between $A_h$ and $B_i$; but this contradicts the item (1), since there is an edge between $A_h$ and $A_i$.
\end{proof}

\begin{lemma}\label{thm:vabbaslemma2}
Let $G$ be a domain. Let $H$ be a single connected component of $G^{0-1}$, and let  $A \cup B$ be the bipartition of the vertices of $H$. There exists a basic $1$-cover $\aa$ of $G$ that yields $1$ on $A$, $0$ on $B$, and such that the cover $\bb$ defined as 
\[b_i \ =\ \left \{
	\begin{array}{cc}
	1 - a_i &\hbox{ if } i \in H, \\
	a_i &\hbox{ otherwise.}\\
        \end{array} 
\right.\]
is a basic $1$-cover that yields $1$ on $B$ and $0$ on $A$.
\end{lemma}

\begin{proof} Let $H_1, \ldots, H_k$ denote the other connected components of  $G^{0-1}$, and let $A_i \cup B_i$ be the bipartition of the vertices of $H_i$. By Lemma \ref{thm:vabbaslemma1} [item (1)], no two points in $B_i$ are adjacent. If in $G$ there is an edge from $A$ to some $A_i$, then $B_i$ is not connected to $B$ by any edge (cf.~Lemma \ref{thm:vabbaslemma1}, item (3)); if in addition there are edges from $A$ to some $A_j$ with $j \neq i$, 
by \ref{thm:vabbaslemma1}, item (5), there is no edge between $B_i$ and $B_j$ either. Therefore, the vector that yields 
\begin{compactitem}
\item $0$ on $B$,
\item $0$ on all the $B_i$'s such that $A_i$ is connected with an edge to $A$, and
\item $1$ everywhere else, 
\end{compactitem}
is a $1$-cover of $G$. If this $1$-cover is basic, call it $\aa$; otherwise, decompose it into the sum of a basic $1$-cover $\aa$ and a $0$-cover $\cc$. In any case, $\aa$ satisfies the desired properties.
\end{proof}

\begin{definition}
By \textit{norm} of a $k$-cover we mean the sum of all its entries. We denote this as
$|a|~:=~\sum_{i=1}^n~a_i$.
\end{definition}

\begin{thm} \label{thm:MSC} Let $G$ be a graph with $n$ vertices, all of them non-isolated. Then the following are equivalent:

\begin{compactenum}[(1)]
\item $G$ satisfies MSC;
\item every basic $k$-cover of $G$ has norm $\frac{kn}{2}$;
\item for any $k$, the norm of all basic $k$-covers of $G$ is the same;
\item $G$ is an unmixed domain;
\item $\AG$ is a domain, and $I(G)$ is unmixed;
\item every connected component of $G^{0-1}$ is a $K_{a,a}$, for some positive integer $a$;
\item $G^{0-1}$ admits a matching;
\item $G$ admits a matching, and all the basic $1$-covers of $G$ have exactly $\frac{n}{2}$ ones.
\end{compactenum}
\end{thm}

\begin{proof}
\begin{compactitem}
\item[(1) $\Rightarrow$ (2):] the matching consists of $\frac{n}{2}$ edges, so if we show that for every edge $\{i, j \}$ of the matching and for every basic $k$-cover one has $a_i + a_j =k$, we are done. But this follows from Lemma \ref{thm:lemma1}, since $\{i',j'\}$ must be an edge of $G$ whenever $\{i,i'\}$ and $\{j,j'\}$ are edges of $G$. 
\item[(2) $\Rightarrow$ (3):] obvious.
\item[(3) $\Rightarrow$ (4):] setting $k=1$ we get the definition of unmixedness. Now, denote by $f(k)$ the norm of any basic $k$-cover. Since twice a basic $1$-cover yields a basic $2$-cover, $2 \cdot f(1) = f(2)$; and in general $k \cdot f(1) = f(k)$. Suppose that an $(s+2t)$-cover $\aa$ can be written as the sum of $s$ basic $1$-covers and $t$ basic $2$-covers. The norm of $\aa$ can be computed via its summands, yielding 
$|\aa| = s \cdot f(1) + t \cdot f(2) = (s+2t) \cdot f(1).$ 
Were $\aa$ non-basic, it could be written as the sum of a basic $(s+2t)$-cover $\bb$ and a $0$-cover $\cc$, whence\\
$ |\aa| = |\bb| + |\cc| \geq  |\bb| + 1 = f(s+2t) + 1 = (s+2t) \cdot f(1) + 1,$
a contradiction. This proves that $G$ is a domain.
\item[(4) $\Leftrightarrow$ (5):] explained in the Introduction. 
\item[(4) $\Rightarrow$ (6):] let $H$ be a connected component of $G^{0-1}$. By Theorem \ref{thm:WSC}, since $G$ is a domain, $H$ is bipartite complete; let $A \cup B$ be the bipartition of the vertices of $H$. Construct the basic $1$-covers $\aa$ and $\bb$ as in Lemma \ref{thm:vabbaslemma2}; they have a different number of ones unless $|A| = |B|$, so by unmixedness we conclude.
\item[(6) $\Rightarrow$ (7), (6) $\Rightarrow$ (8):] obvious.
\item[(7) $\Rightarrow$ (6):] by Lemmas \ref{thm:SC} and \ref{thm:SCG} every connected component of $G^{0-1}$ is either a point or a bipartite complete graph; in order for $G^{0-1}$ to admit a matching, each connected component of $G^{0-1}$ must be of the form $K_{a,a}$, for some integer $a$.
\item[(8) $\Rightarrow$ (1):] let $\mathbb{M}$ be the given matching. In view of Lemma \ref{thm:lemma1}, we only need to show that for every edge $\{i, j \}$ of $\mathbb{M}$ and for every basic $1$-cover $\aa$ one has $a_i + a_j =1$. Yet for any basic $1$-cover $\aa$ one has 
\[ \frac{n}{2} \ = \ \sum_{i=1}^n a_i  \ =\ \sum_{ \{i, j\} \in \mathbb{M}} a_i + a_j, \]
and a sum of $\frac{n}{2}$ positive integers equals $\frac{n}{2}$ only if each summand equals~$1$.
\end{compactitem}
\end{proof}

\begin{corollary}[Villarreal] \label{thm:villarreal}
Let $G$ be a bipartite graph without isolated points. $G$ is unmixed if and only if $G$ satisfies $MSC$.
\end{corollary}

\begin{proof}
In the bipartite case, ``unmixed'' implies ``domain'' (see e.g. \cite{HHT} or \cite{BCV}). So the condition (3) in Theorem \ref{thm:MSC} is equivalent to unmixedness.
\end{proof}

The next result shows how to produce many examples of graphs (not necessarily bipartite) that satisfy the assumptions above.

\begin{thm} \label{thm:technique1}
Let $G$ be an arbitrary graph. 
\begin{compactitem}
 \item Let $G^+$ be the graph obtained by attaching a pendant to each vertex of $G$. Then $G^+$ satifies MSC. Moreover, $G^+$ is bipartite if and only if $G$ is bipartite.
\item Let $G'$ be the graph obtained from $G$ attaching a pendant to each of those vertex of $G$ that appear as isolated vertices in $G^{0-1}$. Then $G'$ satisfies WSC. Moreover, $G'$ satisfies MSC if and only if $G^{0-1}$ is unmixed.
\end{compactitem}
\end{thm}

\begin{proof}
Let us show the second item first. By definition of $G^{0-1}$ (and by Lemma \ref{thm:lemma1}), the isolated vertices of $G^{0-1}$ are exactly the vertices of $G$ at which the WSC property does not hold. By attaching a pendant $j$ at the vertex $i$, the property ``if $\{i,i'\}$ and $\{j,j'\}$ are edges, then $\{i', j'\}$ is also an edge'' holds true trivially, since $j'$ must coincide with $i$. 

Of course, in a matching of $G'$ each pendant should be paired with the vertex it was attached to. By Theorem \ref{thm:MSC}, $(G')$ satisfies MSC if and only if $(G')^{0-1}$ has a matching, if and only if the graph obtained removing all isolated vertices from $G^{0-1}$ has a matching. This happens if and only if each connected component of $G^{0-1}$ is either a single point or a $K_{a,a}$, for some positive integer $a$. This characterizes unmixedness within the class of graphs satifying the SC property.

To prove the first item, label $1^{+}, 2^{+}, \ldots, n^{+}$ the pendants attached to $1, 2, \ldots, n$, respectively. The requested matching is $\{1, 1^{+}\}, \{2, 2^{+}\},\ldots, \{n, n^{+}\}$. The possible presence of an odd-cycle in $G$ reflects in the presence of the same odd cycle in $G'$.
\end{proof}

\par \noindent
\textbf{Examples and Remarks.}\\
\begin{compactenum}
\item The complete graph $G=K_4$ is unmixed and has a matching, but it does not satisfy MSC (it is not even a domain in fact). Every $1$-cover has three ones, while $\frac{n}{2} = 2$. 
Note that the property ``all basic $1$-covers have norm $\frac{n}{2}$'' is strictly stronger than unmixedness, while the property ``for each $k$, all basic $k$-covers have norm $\frac{kn}{2}$'' is equivalent to ``for each $k$, the norm of all basic $k$-covers is the same''. 
\item In Theorem \ref{thm:MSC}, the assumption that no vertex is isolated was introduced only to simplify the notation. In general, an arbitrary graph $G$ (with at least one edge) is an unmixed domain if and only if the reduced graph $G^{red}$, obtained by deleting from $G$ the isolated points, is an unmixed domain. Clearly the basic $1$-covers of $G$ will have exactly $\frac{|G^{red}|}{2}$ ones, and so on.
\item In view of Proposition \ref{thm:technique1} one might think that attaching pendants will make it more likely for a given graph to be a domain. However, let $G$ be the graph with edges $\{1, 2\}$, $\{2, 3\}$, $\{3, 4\}$, $\{4, 1\}$, and $\{2, 5\}$ (a square with a pendant attached). This $G$ is a domain, yet if we attach a pendant to the vertex $3$ the resulting graph is not a domain. 
\item A basic $2$-cover that cannot be the sum of two $1$-covers is said to be \emph{indecomposable}. Bipartite graphs have no indecomposable $2$-covers \cite{HHT}, so in some sense the number of indecomposable $2$-covers of a graph measures its ``distance'' from being bipartite. Suppose $G$ contains an odd cycle and a vertex $i$ none of whose neighbours is part of the cycle. One can see then that $G$ admits a basic $2$-cover $\aa$ that yields $0$ on $i$ and $1$ on the cycle; such an $\aa$ is indecomposable, i.e. it cannot be the sum of two $1$-covers. Now, the property of containing an odd cycle and a ``distant'' vertex is clearly inherited by $G^{+}$, which satisfies MSC. This way one can see that the distance of an unmixed domain from being bipartite can be arbitrarily large.
\end{compactenum}

\bigskip
\par \noindent
\textbf{Acknowledgements.}\\
We would like to thank Aldo Conca and J\"urgen Herzog for the kind support and encouragement.


\begin{thebibliography}{99}
\addcontentsline{toc}{chapter}{Bibliografia}
\bibitem{BCV} B. Benedetti, A. Constantinescu, M. Varbaro, \textit{Dimension, depth and zero-divisors of the algebra of basic $k$-covers of a graph},  Le Matematiche \textbf{LXIII} (2008), pp. 117-156.
\bibitem{HHT} J. Herzog, T. Hibi, N.V. Trung, \textit{Symbolic powers of monomial ideals and vertex cover algebras}, Adv. in Math. \textbf{210} (2007), pp. 304-322.
\bibitem{V} R. H. Villarreal, \textit{Unmixed bipartite graphs}, Revista Colombiana de Matem\'{a}ticas \textbf{41} (2007), pp. 393--395. 
\end{thebibliography}
\end{document}